\documentclass[runningheads]{llncs}
\usepackage{}
\usepackage{graphicx}
\usepackage{float}
\usepackage{amsfonts,amssymb,amsbsy,amsmath}
\usepackage{mathrsfs}
\usepackage{subcaption}
\usepackage{bm}
\usepackage{enumitem}
\begin{document}
%
% \title{Contribution Title\thanks{Supported by organization x.}}
\title{Tailored Finite Point Operator Networks for Interface problems}
\titlerunning{TFPONets for Interface problems}
% % If the paper title is too long for the running head, you can set
% % an abbreviated paper title here
% %
% \author{Anonymous submission}
% \institute{Anonymous submission}
\author{Ye Li\inst{1} \and Ting Du\inst{2}\and  Zhongyi Huang\inst{2}}

\authorrunning{Y. Li et al.}
\institute{Nanjing University of Aeronautics and Astronautics, Nanjing, China \and
Tsinghua University, Beijing, China\\
\email{yeli20@nuaa.edu.cn, dt20@mails.tsinghua.edu.cn, zhongyih@tsinghua.edu.cn}}
\maketitle              % typeset the header of the contribution

\begin{abstract}
Interface problems pose significant challenges due to the discontinuity 
of their solutions, particularly when they involve singular 
perturbations or high-contrast coefficients, resulting in intricate 
singularities that complicate resolution. The increasing adoption 
of deep learning techniques for solving partial differential 
equations has spurred our exploration of these methods for 
addressing interface problems. In this study, we introduce 
Tailored Finite Point Operator Networks (TFPONets) as a novel 
approach for tackling parameterized interface problems. 
Leveraging DeepONets and integrating the Tailored Finite 
Point method (TFPM), TFPONets offer enhanced accuracy in 
reconstructing solutions without the need for intricate equation 
manipulation. Experimental analyses conducted in both one- and two-dimensional scenarios reveal that, in comparison to existing methods such as DeepONet and IONet, TFPONets demonstrate superior learning and generalization capabilities even with limited locations.
\end{abstract}

\section{Introduction}
Elliptic and parabolic problems featuring discontinuous coefficients or sources are prevalent across diverse fields, including fluid mechanics \cite{fadlun2000combined,sussman1999efficient}, materials science \cite{LIU2020109017,WANG2019117}, electromagnetics \cite{hesthaven2003high}, and biomimetics \cite{ji2018finite}. Typically, while solutions to interface problems exhibit smoothness within homogenous material or fluid subdomains, the global solution regularity can be significantly compromised, manifesting singularities at interfaces \cite{babuvska1970finite,huang2009tailored,kellogg1971singularities}. The variable nature of singularities and the complex geometries at interfaces render traditional finite element or finite difference methods less effective for achieving precision. 

The integration of deep learning into computational science has initiated a transformative approach to solving partial differential equations (PDEs). Innovations such as physics-informed neural networks \cite{raissi2019pinn}, the Deep Ritz Method \cite{yu2018deep}, and the Deep Galerkin Method \cite{sirignano2018dgm} leverage neural networks to navigate these challenges. Moreover, operator learning strategies, including the Fourier Neural Operator (FNO) \cite{li2020fourier} and Deep Operator Networks (DeepONets) \cite{lu2021learning}, are emerging as powerful tools capable of approximating mappings between infinite-dimensional spaces. Addressing interface problems—where solutions demonstrate piecewise continuity and are subject to jump conditions at interfaces—poses a challenge for vanilla neural networks methods, often leading to poor performance near the interfaces.

In the manuscript, we introduce Tailored Finite Point Operator Networks (TFPONets), a novel integration of the Tailored Finite Point method (TFPM)\cite{huang2009tailored} with DeepONets, designed to resolve parameterized interface problems with high precision. TFPONets are adept at handling not only general interface problems but also singularly perturbed and high-contrast interface problems, as our experimental results demonstrate robust performance even in these complex scenarios. While our model utilizes DeepONets, its architecture is versatile, allowing for the substitution of other operator networks or neural operators, such as the Fourier Neural Operator (FNO). The key contributions of our work are as follows:  
\begin{itemize}
    \item We introduce TFPONets, a method that combines TFPM for accurate solution reconstruction with the learning capabilities of DeepONets, eliminating the need for complex numerical PDEs discretization and calculation. By incorporating the prior knowledge of basis functions into the neural network design, our method accelerates the learning process without imposing additional training complexity.
    \item TFPONets preserve the meshless nature inherent in DeepONets, allowing for the utilization of either randomly sampled or fixed grid points without constraint. Our models can be trained on coarse grids while delivering accurate predictions for functions beyond the training set and finer locations, exhibiting remarkable generalization capability.
    \item We demonstrate the efficacy of our model across one-dimensional and two-dimensional interface problems, including those with singular perturbations and high-contrast coefficients. Our findings indicate that TFPONets, even with limited locations in the training dataset, can resolve complex interface problems, effectively capturing boundary and interior layers.
\end{itemize}

\section{Preliminaries}
\subsection{Interface Problems}\label{sec:interface problems}

In this manuscript, we develop learning methods to solve an elliptic interface problem specified by the following equation:
\begin{equation}
-\nabla\cdot(a\nabla u)+bu=f,\ \text{in }\Omega/ \Gamma,\ \left[u\right]_{\Gamma}=g_D,\ \left[a\nabla u\cdot\bm{n}\right]_{\Gamma}=g_N,\ u|_{\partial\Omega}=h,
    \label{equ:interface problem ver1}
\end{equation}
where $\Omega$ is the domain comprising $N$ subdomains $\{\Omega_i, i=1,\dots,N\}$ separated by the interface $\Gamma$. Here, $a(\bm{x}) > 0$ and $b(\bm{x}) \geq 0$ are piecewise smooth functions. The terms $[\cdot]{\Gamma}$ and $\left[ a \nabla u \cdot \bm{n} \right]{\Gamma}$ represent the jump across the interface, i.e., for a point $\bm{x}0 \in \Gamma$ between $\Omega_i$ and $\Omega{i+1}$,
\[\left[u\right](\bm{x}_0)=\lim_{\bm{x}\in\Omega_{i+1},\bm{x}\to\bm{x}_0}u(\bm{x})-\lim_{\bm{x}\in\Omega_i,\bm{x}\to\bm{x}_0}u(\bm{x}),\]
\[\left[a\nabla u\cdot\bm{n}\right](\bm{x}_0)=\lim_{\bm{x}\in\Omega_{i+1},\bm{x}\to\bm{x}_0}a(\bm{x})\nabla u(\bm{x})\cdot\bm{n}-\lim_{\bm{x}\in\Omega_i,\bm{x}\to\bm{x}_0}a(\bm{x})\nabla u(\bm{x})\cdot\bm{n}.\]

We also address two more challenging types of interface problems:
\begin{itemize}
    \item \textbf{singularly perturbed interface problem.} When the coefficient function $a(x)$ in the equation is a small positive constant, the resulting solutions display rapid changes in the form of boundary, interior, and corner layers.
    \item \textbf{high-constrast interface problem.} These arise when the coefficient function $a(x)$, while smooth within each subdomain, exhibits significant variations between different subdomains. An example is the modeling of elastic materials with thin layers.
\end{itemize}

\subsection{Current Approaches}

\subsubsection{Numerical solver.} Numerical methods like the immersed boundary (IB) method \cite{peskin1977numerical}, immersed interface method (IIM) \cite{leveque1994immersed}, ghost fluid method (GFM) \cite{liu2000boundary}, and tailored finite point method (TFPM) \cite{huang2009tailored} have been developed to address these problems, each improving the accuracy of capturing interface discontinuities and complex conditions. Nevertheless, general numerical methods necessitate meticulous region discretization and equation formulation.
   
\subsubsection{Neural Network solver.} In the realm of neural networks, piecewise networks have been employed to approximate solutions by assigning a network to each subdomain and using jump conditions as penalty terms in the loss function \cite{guo2021deep,he2022mesh}. Despite their effectiveness, these methods can be complicated by the need to balance various loss terms during optimization. Domain decomposition methods, such as the Schwarz method \cite{lions1988schwarz}, have been adapted with neural networks to ease these challenges \cite{li2019d3m,li2020deep}. Additionally, approaches like the discontinuity capturing shallow neural network (DCSNN) \cite{hu2022discontinuity} and deep Nitsche-type method \cite{wang2020mesh} have been proposed to manage high-contrast discontinuous coefficients and complex boundary conditions. However, these methods are confined to solving a single equation.

\subsubsection{Neural operator.} For parametric PDE problems involving repeated evaluations across similar inputs, operator learning techniques like Fourier Neural Operators (FNO) and Deep Operator Networks (DeepONets) have shown promise. However, research in solving parametric interface problems is limited, with some studies like IONet \cite{wu2023solving} using DeepONets for interface problems across subdomains. Challenges remain in tackling problems with more pronounced singularities, like singularly perturbed and high-contrast interface problems.

\section{Method}
\subsection{Tailored Finite Point Method}\label{subsection:tfpm}

Theoretical findings from \cite{huang2009tailored} indicate that Tailored Finite Point Method (TFPM) not only resolves general interface problems with high precision but also exhibits uniform convergence with respect to the perturbation parameter. This section outlines the Tailored Finite Point method for solving interface problems, as proposed by \cite{huang2009tailored}, covering both one-dimensional and two-dimensional scenarios.

For the one-dimensional case with domain $\Omega = (\alpha, \beta)$, the transformation $y(x) = \int_{\alpha}^x 1/a(\xi)d\xi$ enables the reformulation of Eq.~\eqref{equ:interface problem ver1}. In the two-dimensional case, where the domain is $\Omega = (\alpha_0, \beta_0) \times (\alpha_1, \beta_1)$, a corresponding transformation $\bm{y} = (y_1, y_2) = \left(\int_{\alpha_0}^{x_1} 1/a(\xi_1, x_2)d\xi_1, \int_{\alpha_1}^{x_2} 1/a(x_1, \xi_2)d\xi_2\right)$ is applied to achieve a similar reformulation. The transformed equation takes the form:
\begin{equation}
-\Delta u(\bm{y}) + c(\bm{y})u(\bm{y}) = F(\bm{y}),
\label{equ:interface problem ver2}
\end{equation}
where $c(\bm{y}) = a(\bm{x}(\bm{y}))b(\bm{x}(\bm{y}))$ represents a transformed coefficient and $F(\bm{y}) = a(\bm{x}(\bm{y}))f(\bm{x}(\bm{y}))$ is the transformed source term. Based on this reformulated equation, we will introduce the Tailored Finite Point method (TFPM).

\subsubsection{One-dimensional Case.}
In a one-dimensional setting with a partition $\{y_j\}$, the solution to Eq.~\eqref{equ:interface problem ver2} within each subinterval $(y_{j-1}, y_j)$ can be approximated by:
\begin{equation}
u_h(y) = \alpha_jA_1^j(y) + \beta_jA_2^j(y) + \int_{y_{j-1}}^{y_j} F(s)G_j(y(x),s)ds,
\label{equ:1d_expansion}
\end{equation}
where $G_j$ is the Green's function associated with the subinterval. The functions $A_1(y)=\{A_1^j(y)\}$ and $A_2(y)=\{A_2^j(y)\}$ serve as a pair of local bases, reflecting the distinct characteristics of the various subintervals. They are piecewise defined for each subinterval as follows:
\begin{equation}  
    A_1^j(y),\ A_2^j(y) = \begin{cases}  
        1, \ y,& a_j = b_j = 0, \\  
        \exp(y\sqrt{b_j}),\ \exp(-y\sqrt{b_j}), & a_j = 0,\ b_j \ne 0, \\  
        \text{Ai}(c_h(y)a_j^{-2/3}),\ \text{Bi}(c_h(y)a_j^{-2/3}), & a_j \ne 0,\ b_j \ne 0,  
    \end{cases}  
    \label{equ:1d_basis}  
\end{equation}
where $c_h(y)$ is a piecewise linear approximation of the function $c(y)$, characterized by a slope $a_j$ and an intercept $b_j$ over the $j$-th subinterval. Ai($\cdot$) and Bi($\cdot$) denote the Airy functions of the first and second kind, respectively. The coefficients $\alpha_j$ and $\beta_j$ are determined by a linear system that ensures continuity or enforces jump conditions at the interfaces between subintervals.  

\subsubsection{Two-dimensional Case.}
In two dimensions, we subdivide the domain $\Omega$ into a set of subdomains $\{\Omega^j\}$. By selecting a reference point $\bm{y}^j$ within each $\Omega^j$, the solution to Eq.~\eqref{equ:interface problem ver2} over $\Omega^j$ can be approximated by:
\begin{equation}
u_h(\bm{y})=u_j(r)+\sum\limits_{n=0}^{+\infty}(a_n^jI_n(\mu_j r)\cos n \theta+b_n^jI_n(\mu_j r)\sin n\theta),
    \label{equ:2d expansion}
\end{equation}
where $u_j(r)=F_j\sum\limits_{n=0}^{+\infty}\frac{(\mu_j)^{2n-2}r^{2n}}{4(n!)^2}$. Here, $\mu_j$ is a constant approximation of the function $\sqrt{c(\bm{y})}$ and $F_j$ approximates $F(\bm{y})$ within the $j$-th subdomain. The coordinates $(r, \theta)$ represent the polar coordinates centered at $\bm{y}^j$, and $I_n$ denotes the modified Bessel function of the first kind of order $n$.

Finally, we truncate $u_h$ to include only the leading terms of the series expansion, leaving several unknowns $a_n^j$ and $b_n^j$. These parameters are determined by ensuring $u$ and $\partial u/\partial \bm{n}$ are continuous at $\Omega^j$'s boundary.

\subsection{Deep Operator Networks}
Deep Operator Networks (DeepONets) are designed to learn nonlinear mappings between infinite-dimensional function spaces, particularly adept at solving parametric partial differential equations (PDEs) by establishing relationships between two such spaces. DeepONet operates by processing two inputs: a function $f$ and a location $x$. It consists of two sub-networks: the branch network handles the function $f$ at predetermined points $\{x_i\}$, while the trunk network handles the location $x$. The final output is the dot product of the outputs from these sub-networks, as shown in the Deeponet block in Fig. \ref{fig:tfponet}.

Due to the inherent discontinuity of solutions to interface problems, which exhibit significant variations across different sub-regions, direct application of DeepONet often leads to inaccuracies at the interface. While IONet\cite{wu2023solving} employs multiple DeepONets to handle solution at different sub-regions, it struggles with more singular problems, such as those involving singular perturbations and high-contrast coefficients. To address these challenges, we will propose an enhanced model capable of accurately solving interface problems, including the aforementioned singular cases.

\subsection{Improved Network Structure}

\begin{figure}[ht]  
    \centering
    \includegraphics[width=.7\linewidth]{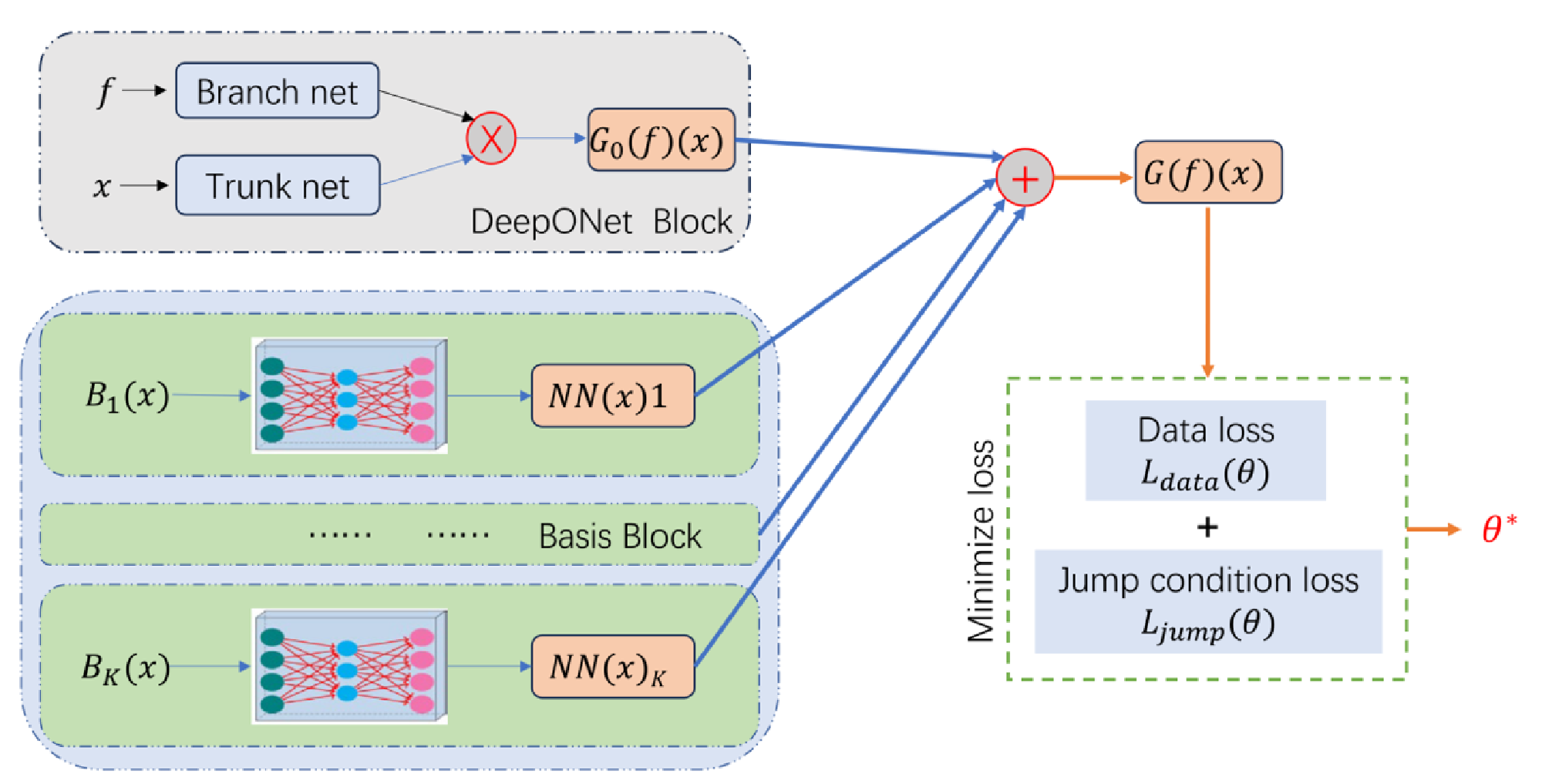}
    \caption{TFPONet integrates a DeepONet with multiple neural networks that receive inputs of the form $(f, x, B_1(x), B_2(x), \dots)$. Specifically, DeepONet is fed with the tuple $(f, x)$, while the $i$-th auxiliary neural network within TFPONet is tasked with processing the local basis $B_i(x)$.}
    \label{fig:tfponet}
\end{figure}

In Section \ref{subsection:tfpm}, we establish that the solution to interface problems can be effectively approximated by a reconstructed solution $u_h$ (see Eq. \eqref{equ:1d_expansion} and \eqref{equ:2d expansion}), which utilizes a series expansion over a finite or infinite set of local basis functions. This approach not only yields high-fidelity approximations but also maintains parameter consistency for singular perturbation problems without necessitating a dense discretization. Consequently, an accurate approximation of $u_h$ directly translates to a precise representation of the interface problem's solution, with the added benefit of a simpler computational framework. This realization underpins the development of the Tailored Finite Point Operator Network (TFPONet), depicted in Fig.~\ref{fig:tfponet}, which augments the vanilla Deep Operator Network (DeepONet).

The TFPONet framework integrates a DeepONet block with multiple basis blocks, each underpinned by a shallow neural network that processes a distinct basis function evaluated at point $x$. The resulting model is encapsulated by the following formulation:
\begin{equation}
\text{TFPONet}(f)(x)=\text{DeepONet}(f)(x)+\sum_{i=1}^K\text{NN}_i(B_i(x)),
\label{equ:tfponet}
\end{equation}
where $\text{DeepONet}(f)(x)$ represents the output from the DeepONet block, and $\text{NN}_i\ (i=1,2,\ldots,K)$ denotes the shallow neural networks corresponding to the ``NNi" blocks in Fig. \ref{fig:tfponet}.

As illustrated in Fig.~\ref{fig:tfponet}, TFPONet comprises a DeepONet block and a series of basis blocks. Each basis block processes a basis function at location $x$, contributing to the approximation of the corresponding term in the reconstructed solution as described by Eq.~\eqref{equ:1d_expansion} or \eqref{equ:2d expansion}. In the one-dimensional case, these terms are $\alpha_jA_1^j(y(x))$ and $\beta_jA_2^j(y(x))$ from Eq.~\eqref{equ:1d_expansion}, while in the two-dimensional case, these terms are $(a_n^j \cos n\theta + b_n^j \sin n\theta)I_n(\mu_j r)$ from Eq.~\eqref{equ:2d expansion}. The DeepONet block is designed to receive the function $f$ and location $x$, enabling the approximation of non-basis terms, exemplified by $\int_{y_{j-1}}^{y_j}F(s)G_j(y(x),s)ds$ in the context of Eq.~\eqref{equ:1d_expansion}, or $u_j(r)$ in Eq.~\eqref{equ:2d expansion}.  This block further serves as an error corrector. Approximation inaccuracies in the basis blocks, stemming from both their intrinsic limitations and the finite number of blocks implementable within TFPONet, necessitate this error correction mechanism. The integration of the DeepONet block harnesses its advanced approximation capabilities to mitigate these errors effectively.

Constructing the local basis functions for basis blocks is straightforward. In the one-dimensional scenario, either a single or a pair of basis blocks suffice, corresponding to inputs $A_1(y(x))$ or $A_1(y(x)),A_2(y(x))$, respectively, as delineated in Eq.~\eqref{equ:1d_basis}. For the two-dimensional case, we employ the $0$-th to $(K-1)$-th order modified Bessel functions of first kind as inputs for $K$ basis blocks. Each basis block is tasked with approximating each local basis term. Consequently, employing a greater number of basis blocks generally enhances the solution's accuracy. However, this also leads to an expanded set of neural networks, which, while potentially improving training accuracy and simplifying the training process, necessitates consideration of the inherent trade-off due to the resultant increase in parameter count.
 
\subsection{Optimization Problem and Loss Function}

For interface problems where the solution exhibits a non-zero jump, indicating discontinuity at the interface, the solution varies significantly across sub-regions. To preserve this discontinuity, a suitable approach is the deployment of distinct, smaller TFPONet for each subdomain, with each network dedicated to its respective region. This composite model is henceforth referred to as TFPONets.

If we partition the specified region $\Omega$ into $N$ sub-regions using the interface $\Gamma$, such that $\Omega = \mathop{\cup}\limits_{i=1}^N\Omega_i$, we address the underlying problem with a model where 
\begin{equation}
   \mathcal{N}(f)(x)=\mathcal{N}_i(f)(x),\text{ if }x\ \in\ \Omega_i,
   \label{equ:model}
\end{equation}
with $\mathcal{N}_i$ being the TFPONet for sub-region $\Omega_i$. The overall approach involves training all subproblems simultaneously, with each subproblem optimizing the parameters $\boldsymbol{\theta}_i$ of $\mathcal{N}_i$ within its region $\Omega_i$. This is achieved by utilizing a total loss function that encompasses both individual subproblem losses and the interface loss.

The model's parameters across TFPONets are $\boldsymbol{\theta}=\mathop{\cup}\limits_{i=1}^N\boldsymbol{\theta}_i$, optimized by minimizing the loss function 
\begin{equation}
    \mathcal{L}(\boldsymbol{\theta})=\mathcal{L}_{data}(\boldsymbol{\theta})+\gamma\cdot\mathcal{L}_{jump}(\boldsymbol{\theta}).
\end{equation}
Here, $\mathcal{L}_{data}$ is the data fit term, computed as
\begin{equation}
    \mathcal{L}_{data}(\boldsymbol{\theta})=\frac{1}{MJ}\sum\limits_{m=1}^{M}\sum\limits_{j=1}^{J}|u_m(x_j)-\mathscr{N}_{\boldsymbol{\theta}}(F_m)(x_j)|^2,
\end{equation}
where the ground truth, $u_m$, is obtained through high-precision numerical methods, specifically using TFPM on a uniform mesh. The samples $F_1,\dots,F_M$ are independently drawn from Gaussian random fields, and the points $x_1,\dots,x_J$ are positioned on a uniform grid within the region. $\mathcal{L}_{jump}$ is a soft constraint enforcing jump conditions at interfaces, with $\gamma$ as the penalizing parameter, as specified by
\begin{equation}
\begin{aligned}
    \mathcal{L}_{jump}(\boldsymbol{\theta})=&\frac{1}{MJ_0}\left(\sum\limits_{m=1}^{M}\sum\limits_{i=1}^{J_0}|[\mathcal{N}_{\boldsymbol{\theta}}(F_m)](x_i^{\gamma})-g_D(x^{\gamma}_i)|^2\right.\\
   &\left.+\sum\limits_{m=1}^{M}\sum\limits_{i=1}^{J_0}|\left[a\nabla \mathcal{N}_{\boldsymbol{\theta}}(F_m)\cdot\bm{n}\right](x^{\gamma}_i)-g_N(x^{\gamma}_i)|^2\right),
\end{aligned}
\end{equation}
where $[\cdot](x^{\gamma}_i)$ denotes the jump across the interface at $x^{\gamma}_i$, similar to the definition in section \ref{sec:interface problems}.

\section{Experiments}

\subsubsection{Experiment details}. 
In this section, we explore the application of TFPONets to various interface problems encompassing both one-dimensional and two-dimensional scenarios, including those with singular perturbations and high-contrast coefficients. Additionally, we conduct a comparative analysis of our experimental results with those obtained using IONet \cite{wu2023solving}, which is an improvement based on the vanilla DeepONet. IONet employs a DeepONet in each sub-region to handle the solution's singularity at the interface. If the solution is continuous across the entire region, this method degenerates to the vanilla DeepONet. The training dataset is composed of numerous triplets $(f,x,u)$, with $f$ samples being independently drawn from Gaussian random fields to serve as the input function. Here, $x$ denotes a location, chosen without any specific constraints, and for subsequent experiments, uniform grid points are utilized as these points. The solution $u$ of the governing equations at each location $x$ for a given input function $f$ is derived using TFPM. When indicating that the training dataset comprises $n_1 \times n_2$ triplets $(f,x,u)$, we mean that the input data includes $n_1$ random samples $f$ and $n_2$ locations $x$, resulting in a corresponding output dataset with $n_1 \times n_2$ $u$.

\subsubsection{TFPONets Settings}. In all experiments, each TFPONet within the respective sub-region utilizes two basis blocks. Specifically, for one-dimensional scenarios, $A_1(y(x))$ and $A_2(y(x))$—as defined by Eq.~\eqref{equ:1d_basis}—serve as the inputs for these blocks, while for two-dimensional scenarios, the inputs are the zeroth and first-order modified Bessel functions of the first kind. In the one-dimensional scenario, both the branch net, trunk net, and basis block within the DeepONet block utilize a fully connected neural network with a single hidden layer. Additionally, in the two-dimensional setting, they all employ neural networks with two hidden layers. The activation functions across all layers are ReLU.

\subsubsection{One-dimensional Settings}. Our initial experiments address the following one-dimensional problem:
\begin{equation}
-(au')'+bu=f,\ \text{in }\Omega/ \Gamma,\ \left[u\right]_{\Gamma}=g_D,\ \left[u'\right]_{\Gamma}=g_N,\ u|_{\partial\Omega}=0.
    \label{equ:1d_example}
\end{equation}
where $\Omega=(0,1)$ and the interface, denoted by $\Gamma = \{0.5\}$, consists of a single point. Unless otherwise indicated, our goal is to infer the mapping $f \mapsto u$ using multiple sets of input-output data pairs.

\subsubsection{Example 1}. We explore a straightforward scenario wherein both the solution to the problem and its derivative demonstrate continuity at the interface $\Gamma$, that is, $g_D = g_N = 0$. The coefficient $b(x)$ and the source term $f(x)$ in Eq.~\eqref{equ:1d_example} are piecewise functions given by:
    \begin{equation*}
        b(x)=\left\{\begin{array}{lr}
        2x+1,&  x\in [0,0.5),\\
        2(1-x)+1,& x\in [0.5, 1],
    \end{array}\right.\quad
    f(x)=\left\{\begin{array}{lr}
        f_1(x),&  x\in [0,0.5),\\
        f_2(x),& x\in [0.5, 1],
    \end{array}\right.
\end{equation*}
where $f_1(x)$ and $f_2(x)$ are functions drawn from two specified distribution spaces. The discontinuity in these functions or their derivatives makes the solution appear singular at the interface, even though the solution $u$ and its derivatives are continuous. Our objective is to deduce the mapping from $f$ to $u$.

We assess the efficacy of our proposed TFPONets on problems characterized by singular perturbations and high-contrast coefficients. Given the continuity of the solution across the interface, a single TFPONet suffices for the entire domain. This approach will be benchmarked against the vanilla DeepONet.

\begin{figure}[!tbh]
\centering
\begin{minipage}{0.325\textwidth}
\centering
\includegraphics[width=\textwidth]{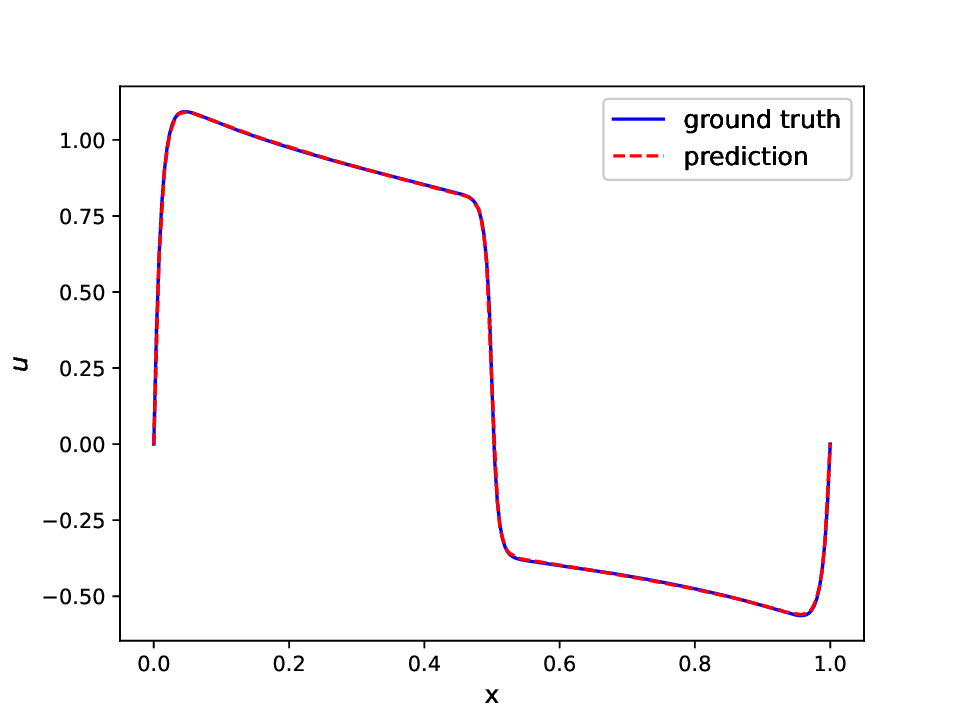}\\
\scriptsize{a)}
\end{minipage}
\begin{minipage}{0.325\textwidth}
\centering
\begin{minipage}{\textwidth}
\centering
\includegraphics[width=\textwidth]{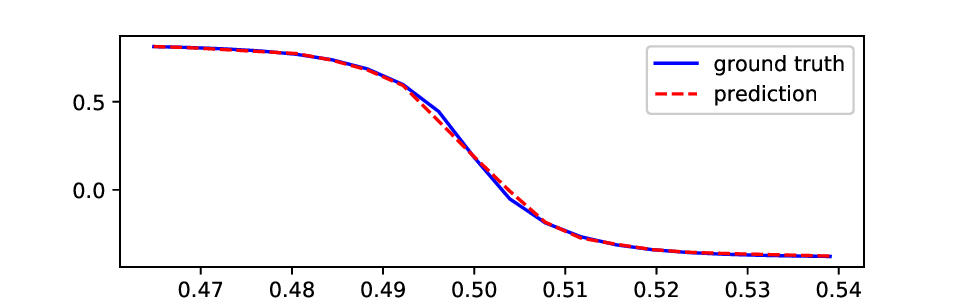}\\
\scriptsize{b)}
\end{minipage}
\\
\begin{minipage}{\textwidth}
\centering
\includegraphics[width=\textwidth]{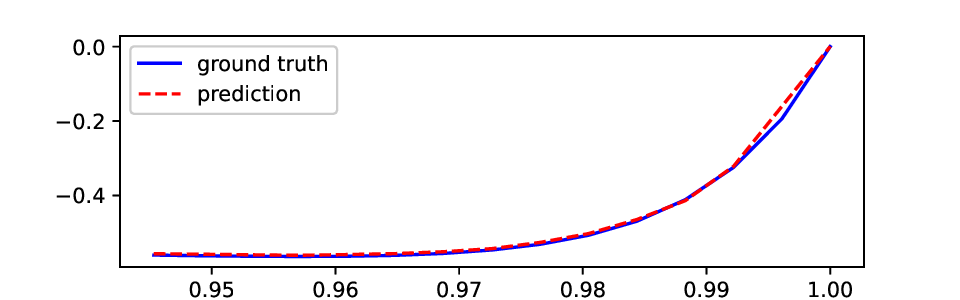}\\
\scriptsize{c)}
\end{minipage}
\end{minipage}
\begin{minipage}{0.315\textwidth}
\centering
\includegraphics[width=\textwidth]{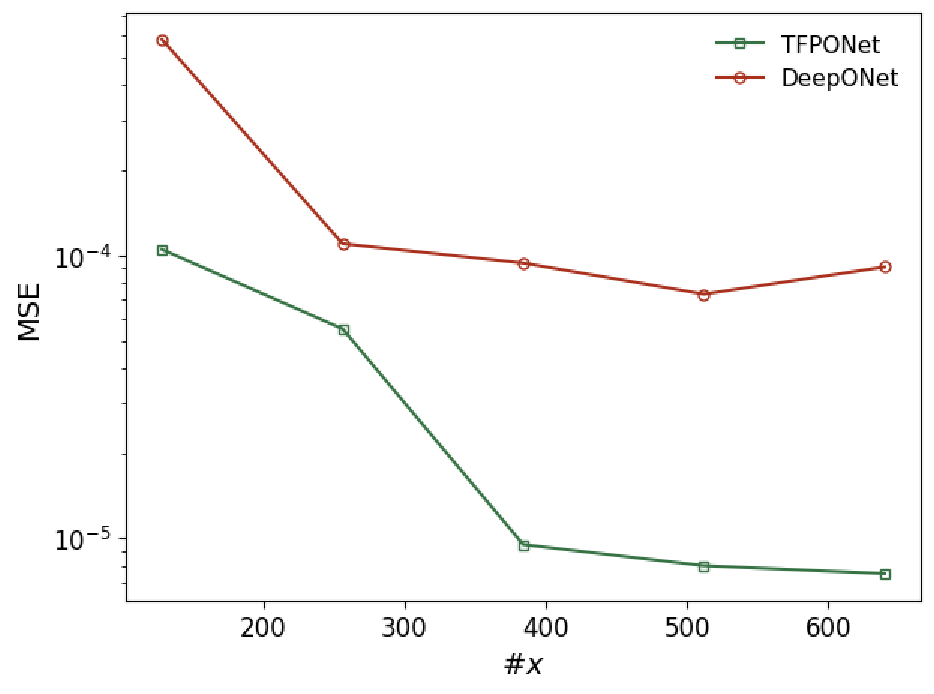}\\
\scriptsize{d)}
\end{minipage}
\caption{a): Ground truth and TFPONet's prediction for Eq.~\eqref{equ:1d_example} with $a(x)=0.0001$ on a random $f$ and 257 locations. b), c): Zoomed-in view. d): Test error (MSE) vs. Training resolution $(\# x)$ for TFPOnet and DeepONet.}
\label{Fig:ex1}
\end{figure}

\begin{figure}[!tbh]
\centering
\begin{minipage}{0.325\textwidth}
\centering
\includegraphics[width=\textwidth]{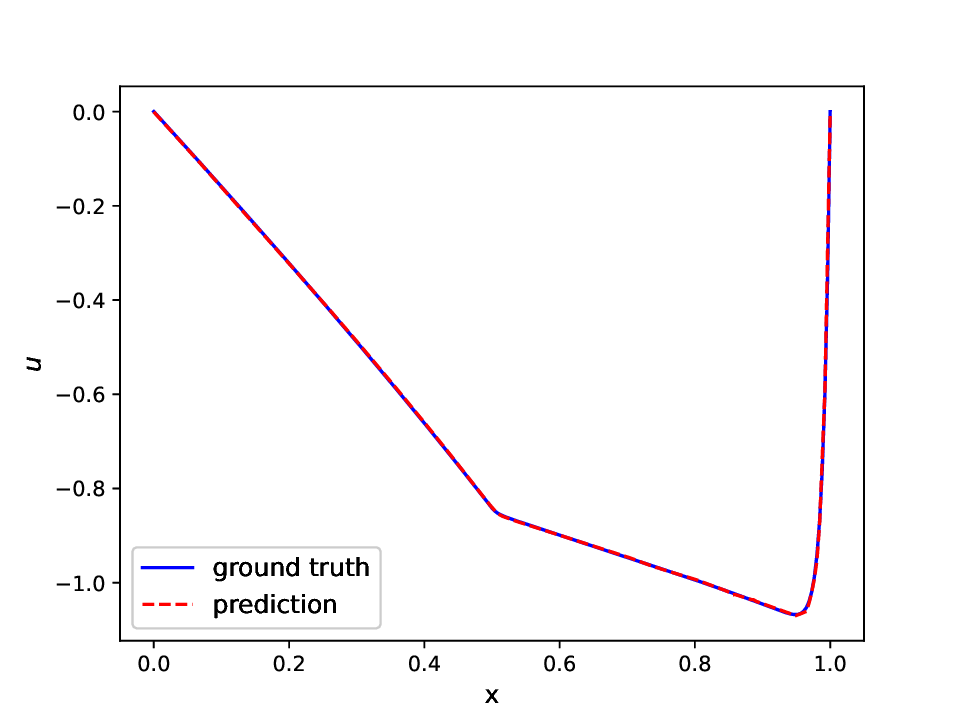}\\
\scriptsize{a)}
\end{minipage}
\begin{minipage}{0.325\textwidth}
\centering
\begin{minipage}{\textwidth}
\centering
\includegraphics[width=\textwidth]{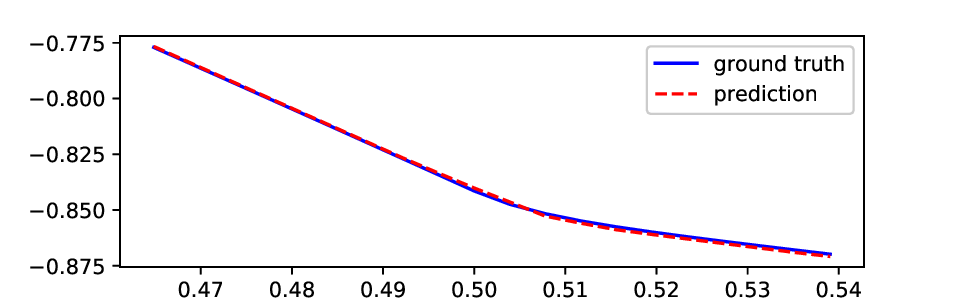}\\
\scriptsize{b)}
\end{minipage}
\\
\begin{minipage}{\textwidth}
\centering
\includegraphics[width=\textwidth]{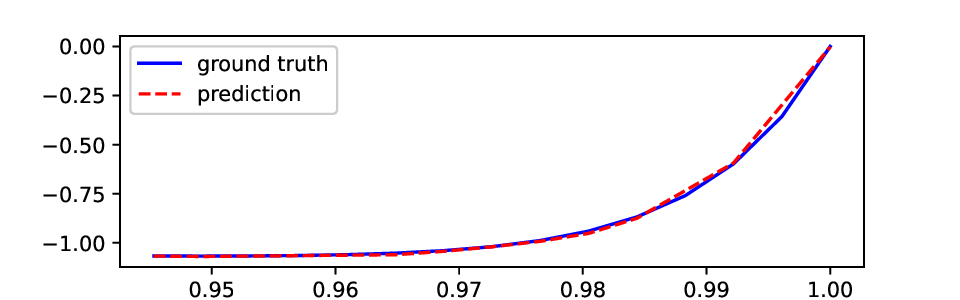}\\
\scriptsize{c)}
\end{minipage}
\end{minipage}
\begin{minipage}{0.315\textwidth}
\centering
\includegraphics[width=\textwidth]{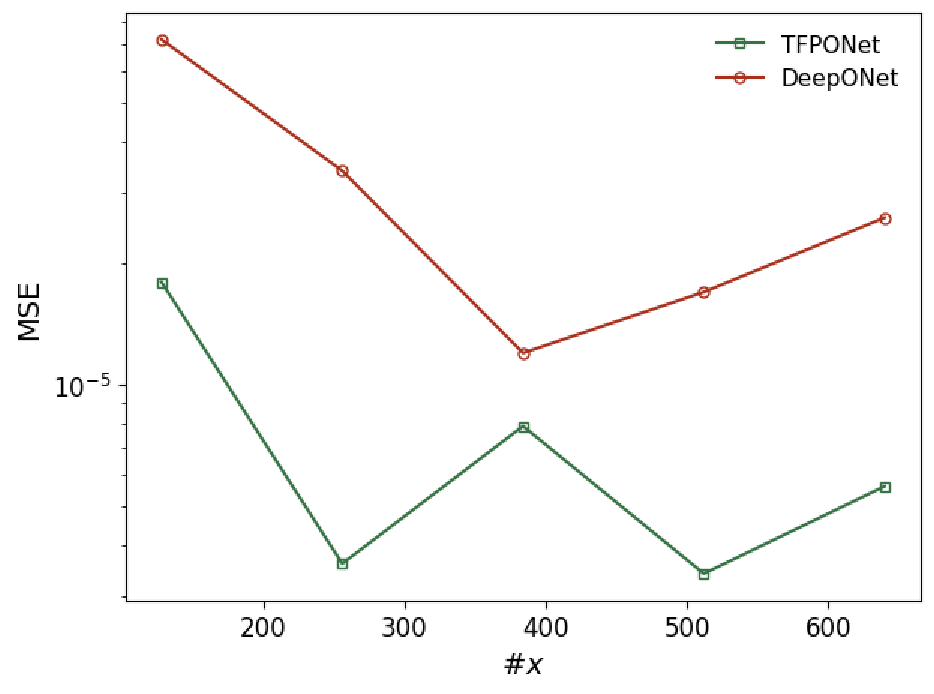}\\
\scriptsize{d)}
\end{minipage}
\caption{a):  Ground truth and TFPONet's prediction for Eq.~\eqref{equ:1d_example}, where $a(x) = 1$ for $x \in [0, 0.5)$ and $a(x) = 0.0001$ for $x \in [0.5, 1]$, across 257 locations using a random $f$. b), c): Zoomed-in view. d): Test error (MSE) vs. Training resolution $(\# x)$ for TFPOnet and DeepONet.}
\label{Fig:ex1_2}
\end{figure}

In the singular perturbation case, setting the coefficient $a(x)$ in Eq.\eqref{equ:1d_example} to 0.0001 induces rapid solution transitions within interior and boundary layers, as illustrated in Fig.\ref{Fig:ex1}$a)$. For high-contrast problems, we define $a(x) = 1$ for $x \in [0, 0.5)$ and $a(x) = 0.0001$ for $x \in [0.5, 1]$, resulting in solutions with distinct singularities at $x = 0.5$ and $x = 1$, shown in Fig.\ref{Fig:ex1_2}$a)$.  

We train the TFPONet using datasets comprising $1000 \times 129$ triples $(f, x, u)$ for both problems. To assess performance, we test the trained model at 257 locations using a random $f$ sample that is independent from the training data. Subplots a), b), and c) in Figs.~\ref{Fig:ex1} and \ref{Fig:ex1_2} compare the predictions with the ground truth, highlighting the model's accuracy in critical regions. Then we train both TFPONet and DeepONet across varying position counts—129, 257, 385, 513, and 641—and assessed them on 100 distinct $f$ samples at 1001 positions. Using mean square error (MSE) as the evaluation metric, we charted test errors relative to the number of training positions ($\# x$) in Figs.~\ref{Fig:ex1}d) and \ref{Fig:ex1_2}d). The results consistently demonstrate TFPONet's superior learning and generalization over vanilla DeepONet for these interface problems.

\subsubsection{Example 2}. \label{example:example2}We next examine an interface problem as defined by Eq.~\eqref{equ:1d_example}, with the solution exhibiting a discontinuity at the interface $\Gamma=\{0.5\}$, specifically where $g_D=g_N=1$. The coefficients are given by $a(x)=1$ and a piecewise function for $b(x)$:
\begin{equation*}
    b(x) = 5000\text{ for }x\ \in [0, 0.5)\text{ and }b(x) = 100(4+32x)\text{ for }x\ \in (0.5, 1].
\end{equation*}
The source term $f$ is continuous, sampled from a Gaussian random field.

\begin{figure}[ht]  
    \centering  
    \begin{minipage}{.5\textwidth}  
        \centering  
        \includegraphics[width=.95\linewidth]{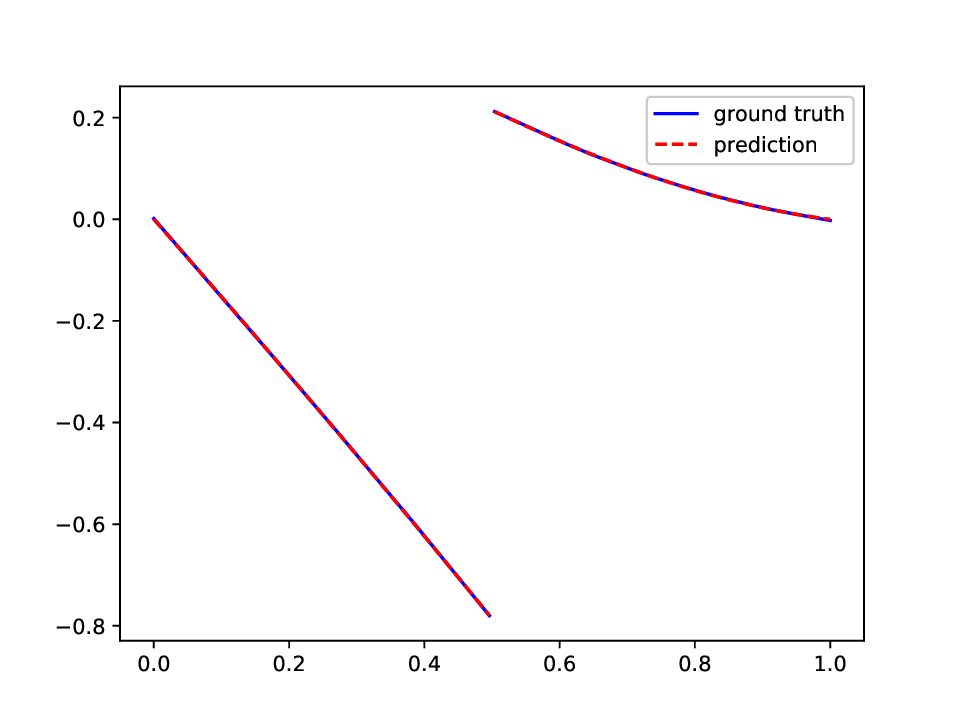}
        \scriptsize{a)}
    \end{minipage}%  
    \begin{minipage}{.5\textwidth}  
        \centering  
        \includegraphics[width=.95\linewidth]{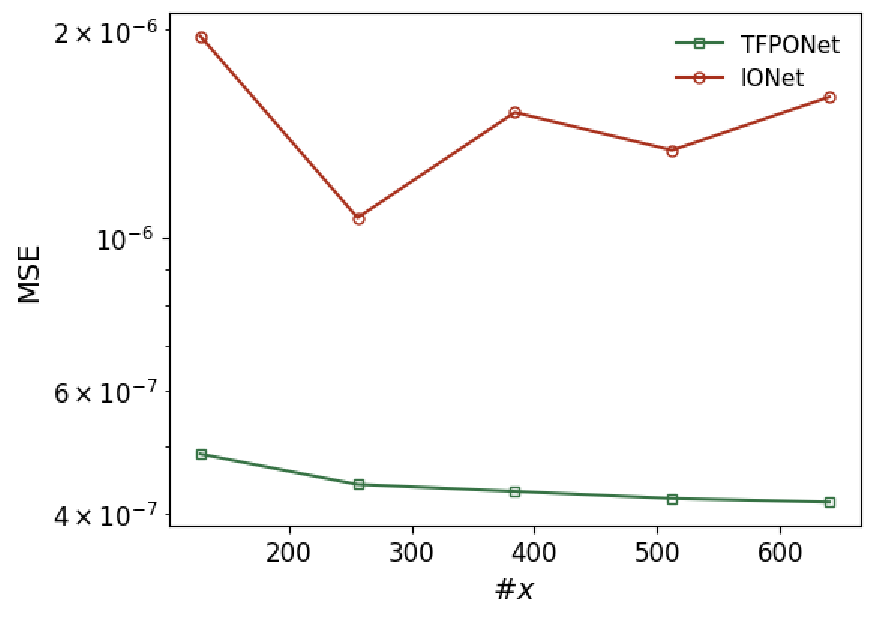}  
        \scriptsize{b)}
    \end{minipage}
    \caption{a): Ground truth and TFPONet's prediction for Example 2. b): Test error (MSE) vs. Training resolution $(\# x)$ for TFPONets and IONet.}
    \label{Fig:ex2}
\end{figure}

To train the TFPONets, which comprise two TFPONets each dedicated to one of the subdomains $[0, 0.5)$ and $(0.5, 1]$, we utilized a dataset consisting of $1000 \times 128$ triplets $(f, x, u)$, with 64 locations uniformly distributed within each subdomain. The trained network is then evaluated using a random $f$ sample at 256 locations. The results, displayed in subplot a) of Fig. \ref{Fig:ex2}, illustrate the prediction versus the ground truth.

We proceed to train both TFPONets and IONet across varying position counts: 128, 256, 384, 512, and 640. Subsequently, we evaluate their performance on 100 distinct $f$ samples at 1001 positions. Subplot b) of Fig. \ref{Fig:ex2} illustrates the testing errors relative to the training resolution. This analysis showcases the robust generalization capabilities of our model, coupled with its remarkable accuracy in handling general interface problems.

\subsubsection{Two-dimensional Settings}. We consider the following equation:
\begin{equation}
    -a\Delta u+bu=0,\ \bm{x}\in\Omega/ \Gamma,\ \left[u\right]_{\Gamma}=1,\ \left[\nabla u\cdot\bm{n}\right]_{\Gamma}=0,\ u|_{\Gamma_1}=0,\ u|_{\Gamma_2\cup\Gamma_3}=f.
    \label{equ:2d_example}
\end{equation}
Here $\Omega=(-1,1)\times(-1,1)$, $a(\bm{x})=0.001$ and $b(\bm{x})$ and $f(x_1)$ are piecewise-defined functions:
\[b(\bm{x})=\left\{\begin{array}{lr}
    (1-x_1^2)^2, &\ x_1<0, \\
    0.001, & \ x_1>0,
\end{array}\right.\quad f(x_1)=\left\{\begin{array}{lr}
    f_1(x_1), &\  x_1<0, \\
    f_2(x_1), &\ x_1>0.
\end{array}\right.\]
The equation is governed by jump conditions on $\Gamma$ and boundary conditions on $\Gamma_1$, $\Gamma_2$, and $\Gamma_3$, where: 
\begin{equation*}
\begin{aligned}
    \Gamma=\{\bm{x}\in\mathbb{R}^2|x_1=0,-1\le x_2\le1\},\ \Gamma_1=\{\bm{x}\in\mathbb{R}^2|x_1=\pm1,-1\le x_2\le1\},\\
    \Gamma_2=\{\bm{x}\in\mathbb{R}^2|-1\le x_1<0,x_2=\pm1\},\ \Gamma_3=\{\bm{x}\in\mathbb{R}^2|0< x_1\le1,x_2=\pm1\}.
\end{aligned}
\end{equation*}
\subsubsection{Example 3}
\begin{figure}[!tbh]
\centering
\begin{minipage}{0.34\textwidth}
\centering
\includegraphics[width=\textwidth]{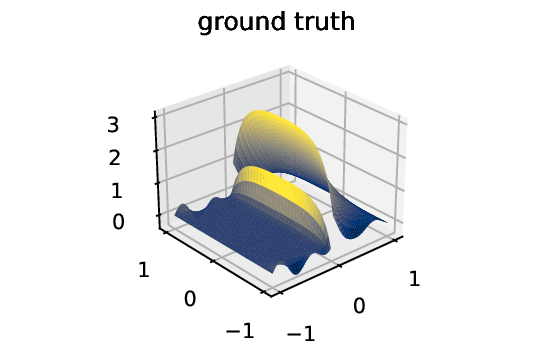}\\
\scriptsize{a)}
\end{minipage}
\begin{minipage}{0.20\textwidth}
\centering
\includegraphics[width=\textwidth]{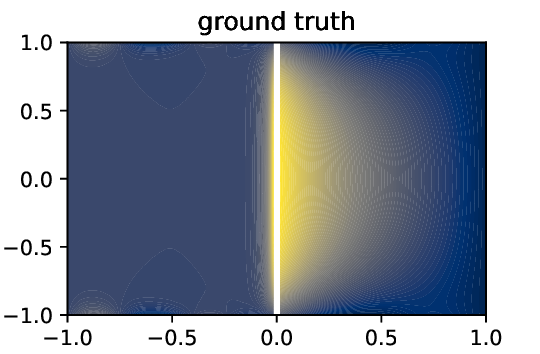}\\
\scriptsize{b)}
\end{minipage}
\begin{minipage}{0.20\textwidth}
\centering
\begin{minipage}{\textwidth}
\centering
\includegraphics[width=\textwidth]{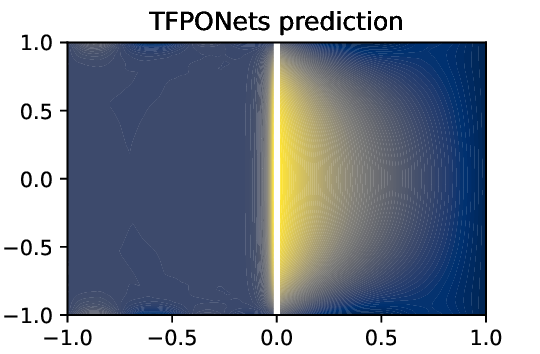}\\
\end{minipage}
\\
\begin{minipage}{\textwidth}
\centering
\includegraphics[width=\textwidth]{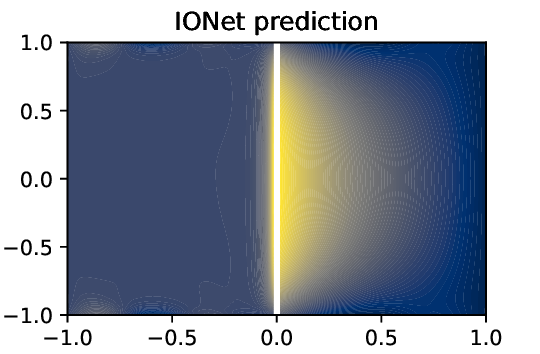}\\
\scriptsize{c)}
\end{minipage}
\end{minipage}
\begin{minipage}{0.21\textwidth}
\centering
\begin{minipage}{\textwidth}
\centering
\includegraphics[width=\textwidth]{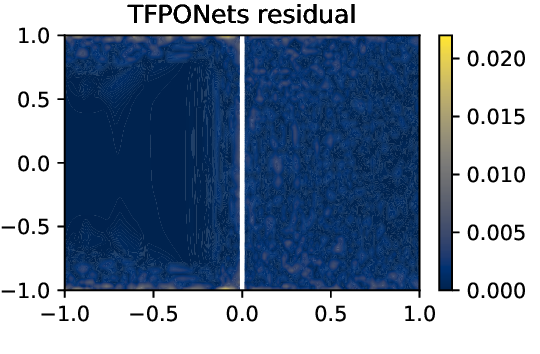}\\
\end{minipage}
\\
\begin{minipage}{\textwidth}
\centering
\includegraphics[width=\textwidth]{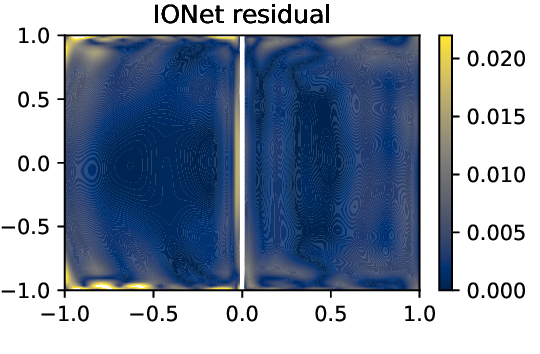}\\
\scriptsize{d)}
\end{minipage}
\end{minipage}
\caption{Performance of TFPONets and IONet on Eq.\eqref{equ:2d_example}. a), b): Ground truth. c): Predictions of TFPONets (above) and IONet (below). d): Absolute Point-wise error of TFPONets (above) and IONet (below).}
\label{Fig:ex3}
\end{figure}

. 
We focus on Eq. \eqref{equ:2d_example} as our research subject, aiming to learn the mapping $f\mapsto u$, where $f_1$ and $f_2$ are independently drawn from two given distribution spaces.
 We randomly select 500 $f$ samples, along with the grid points on a $65\times65$ uniform grid (excluding points on $x_1=0$) as locations. The center point in each generated cell serves as the reference point. Additionally, we employ IONet, utilizing two DeepONets for each subdomain, trained on the same set of $f$ samples and locations. The prediction capabilities of these models for a random $f$ at identical locations are illustrated in Fig.~\ref{Fig:ex3}. Experiments demonstrate that our model excels in capturing the inner layer of solutions to two-dimensional problems compared to IONet.

\section{Conclusion}

In this work, we introduce Tailored Finite Point Operator Networks (TFPONets), a novel integration of the tailored finite point method (TFPM) with Deep Operator Networks (DeepONets), designed to tackle parameterized interface problems. The motivation behind TFPONets stems from the strengths of TFPM, which avoids the need for specialized spatial discretization and complex equation manipulation, using local basis expansion to approximate solutions to interface problems. TFPM has demonstrated high accuracy for standard interface problems and maintains consistency under parameter perturbations in singular perturbation scenarios. By combining TFPM's robustness with DeepONets' powerful learning capabilities, TFPONets effectively address one-dimensional and two-dimensional interface problems, including those with singular perturbations and high-contrast coefficients. Our analysis shows that TFPONets surpass the performance of vanilla DeepONets and IONet in terms of prediction accuracy and generalization, even with spatially sparse training samples.

Despite our experiments being limited to one- and two-dimensional scenarios, the inherent adaptability of our model inspires confidence in its potential applicability to more complex, higher-dimensional problems.

\bibliographystyle{splncs04}
\bibliography{reference}
\end{document}